
\input amstex.tex
\documentstyle{amsppt}
\input epsf

 \advance\hsize4cm
 \advance\vsize3cm
 \advance\voffset-0.5cm
 \advance\hoffset-0.5cm

\define\R{\Bbb R}
\define\Z{\Bbb Z}

\define\emb{\hookrightarrow}
\define\delet{\mathaccent"7017 }
\define\St{\operatorname{St}}
\define\Lk{\operatorname{Lk}}

\define\lk{\operatorname{lk}}

\rightheadtext{Embedding products of graphs...} \NoBlackBoxes

\topmatter
\title Embedding products of graphs
       into Euclidean spaces
\endtitle
\author Mikhail Skopenkov \endauthor
\thanks
This is an improved version of the paper published in Fund. Math. 179 (2003), 191--197.
The author was supported in part by INTAS grant 06-1000014-6277,
Russian Foundation of Basic Research grants 05-01-00993-a,
06-01-72551-NCNIL-a, 07-01-00648-a, President of the Russian
Federation grant NSh-4578.2006.1, Agency for Education and Science grant RNP-2.1.1.7988, and Moebius Contest Foundation for Young Scientists.
\endthanks
\abstract For any collection of graphs $G_1,\dots,G_N$, we find the
minimal dimension $d$ such that the product $G_1\times\dots\times
G_N$ is embeddable into $\R^d$. In
particular, we prove that $(K_5)^n$ and $(K_{3,3})^n$ are not
embeddable into $\R^{2n}$, where $K_5$ and $K_{3,3}$ are the
Kuratowski graphs. This is a solution to a problem by Menger from
1929. The idea of the proof is a reduction to a problem from
so-called Ramsey link theory: we show that any embedding
$L\to S^{2n-1}$, where $L$ is the join of $n$ copies of a $4$-point set, has a pair of linked $(n-1)$-dimensional spheres.
\endabstract
\address
Faculty of Mechanics and Mathematics,
Moscow State University, Moscow, 119992, Russia
\endaddress
\email mikhail$\cdot$skopenkov\,\@\,gmail$\cdot$com \endemail
\endtopmatter

\document

{\bf Introduction.} Our main result is a solution to the Menger
problem from his paper \cite{Men29}:
$$(K_5)^N\not\emb\R^{2N} \qquad\text{and}\qquad (K_{3,3})^N\not\emb\R^{2N}.$$
Hereafter we denote by $K_n$ {\it a complete graph} on $n$
vertices and by $K_{n,n}$ {\it a complete bipartite graph} with $n$
vertices in each part. We write $K\emb L$ if a compact polyhedron $K$ is {\it piecewise linearly embeddable} into a polyhedron~$L$, i.e., there is an injective piecewise linear map $K\to L$. Space $\R^n$ is viewed as a non-compact polyhedron. We get the same result on topological non-embeddability.

Menger posed this problem explicitly for $N=2$  but the title of his paper \cite{Men29} suggests that he was interested in arbitrary $N$ as well. We also solve a more general problem posed by Dranishnikov \cite{Gal92}: given a collection of graphs $G_1,\dots,G_N$, find
the minimal dimension $d$ such that $G_1\times\dots\times
G_N\emb\R^d$.

The topological problem of embeddability is an important one
(e.~g., see \cite {Sch84, ReSk99, ARS01, Sko07}). Our special case
of the problem is interesting because the complete answer can be
obtained and is stated easily, but the proof is nontrivial and
contains interesting ideas. For applications of the result, see the works by Gromov \cite{Gr10} and Lindenstrauss--Tsukamoto \cite{LT14}.

\proclaim{Theorem~1} Let $G_1,\dots, G_n$ be finite connected graphs,
not homemorphic to a point, a closed interval~$I$, and a circle $S^1$. The minimal dimension $d$ such that
$$G_1\times\dots\times G_n\times(S^1)^s\times I^i\emb\R^d$$
equals
$$
d=\cases 2n+s+i, \quad \text{if $i\ne0$ or one of the graphs $G_1,\dots,G_n$ is
planar }
\text {{\rm (}i.e., $\exists k : K_5,K_{3,3}\not\hookrightarrow G_k${\rm),} }\quad &(1)\\
2n+s+1, \quad \text{otherwise. } &(2)
\endcases
$$
\endproclaim

The Menger conjecture is the particular case of the theorem when $G_1=\dots=G_n=K_5$ or $K_{3,3}$.

Hereafter a {\it $d$-dimensional sphere} $S^d$ is any polyhedron that is piecewise linearly homeomorphic to the boundary of a $(d+1)$-dimensional simplex (a slight deviation from the common notation), and $I=[0;1]$.

Theorem~1 remains true in topological category, i.e., with  $K\emb L$ replaced by $K\emb_{TOP} L$.
We write $K\emb_{TOP} L$ if $K$ is compact and there is an injective continuous map $K\to L$.
We first prove Theorem~1
in piecewise linear category and then deduce the topological version from the piecewise linear one.
From now and till that moment we work in the piecewise linear category. See \cite{RS82} for an introduction to the latter.

Theorem~1 was stated (without proof) in \cite{Gal93}, cf. \cite{Gal92}.
The proof of embeddability is trivial (see below).
The non-embeddability has been proved earlier in some particular
cases. For example, it was known that $Y^n\not\emb\R^{2n-1}$,
where $Y$ is a {\it triod} (a graph homeomorphic to the letter ''Y''). A nice proof of this
folklore result is presented in \cite{Sko07}, cf. \cite {ReSk01}.
Also it was known that $K_5\times S^1\not\emb\R^3$ (Tom Tucker,
private communication). In \cite{Um78} Ummel proved that
$K_5\times K_5\not\emb\R^4$ and $K_{3,3}\times K_{3,3}\not\emb\R^4$, thus solving the problem explicitly posed by Menger in \cite{Men29}. That proof contained about 10 pages of calculations involving spectral sequences. We obtain a shorter geometric proof of this result (see Example~2 and Lemma~2 below).
The proof of the non-embeddability in case (2), namely, Lemma~2,
is the main point of Theorem~1 (while case (1) is reduced easily
to a result of van Kampen.)

Our proof of Theorem~1 is quite elementary, in particular, we do
not use any abstract algebraic topology. A popular-science introduction to the method is given in \cite{Sko24}. We use a reduction to a
problem from so-called {\it Ramsey link theory} \cite{S81, CG83,
SeSp92, RST93, RST95, LS98, Neg98, SSS98, T00, ShTa}. The
classical Conway--Gordon--Sachs theorem of Ramsey link theory
asserts that any embedding of $K_6$ into $\R^3$ has a pair of
(homologically) linked cycles. In other words, {\it $K_6$ is not
linklessly embeddable into $\R^3$}. The graph $K_{4,4}$ has the
same property (the Sachs theorem, proved in \cite{S81}). 
Denote by
$\sigma^m_n$ the $m$-skeleton of an $n$-simplex. For a polyhedron
$\sigma$, let $\sigma^{*n}$ be the join of $n$ copies of $\sigma$.
In our proof of Theorem~1, we use the following higher dimensional
generalization of the Sachs theorem.

\proclaim{Lemma 1} Any embedding $(\sigma^0_3)^{*n}\to S^{2n-1}$ has  a pair
of (homologically) linked $(n-1)$-dimensional spheres.
\endproclaim

Lemma~1 follows from Lemma~$1'$ below. Higher dimensional generalizations of the Conway--Gordon--Sachs theorem are known in arbitrary codimension \cite{SeSp92, SSS98, T00}. An open question: is there an $n$-dimensional polyhedron such that any its embedding into $\R^{n+2}$ contains a knotted $n$-dimensional sphere, for $n>1$?

\bigskip

{\bf The easy part of Theorem~1 and some heuristic
considerations.} First, let us prove all the assertions of Theorem~1
except for the nonembeddabili\-ty in case (2).

\demo{Proof of the embeddability in Theorem 1} We need the
following two simple results:

(*) If a compact polyhedron $K\emb\R^d$ and $d>0$, then $K\times I$,
$K\times S^1\emb\R^{d+1}$ (it is sufficient to prove this for $K=I^d$,
for which this is trivial).

(**) For any compact $d$-dimensional polyhedron $K$, the cylinder $K\times
I\emb\R^{2d+1}$ \cite{RSS95}.

Denote $G=G_1\times\dots\times G_n$. By general position,
$G\emb\R^{2n+1}$. If $i\ne0$, then by (**)  $G\times
I\emb\R^{2n+1}$. If, say, $G_1$ is planar, that is, $G_1\hookrightarrow I^2$, then by (**) and (*) we get
$I^2\times G_{2}\times\dots\times G_{n}\emb\R^{2n}$, whence
$G\emb\R^{2n}$. Applying (*) repeatedly, we get the
embeddability assertion in all cases considered. \qed\enddemo

\demo{Proof of the non-embeddability in Theorem 1 in case (1)} Note that any
connected graph, not homemorphic to a point, $I$, and $S^1$, contains a
triod $Y$. So it suffices to prove that $Y^n\times
I^{s+i}\not\emb\R^{2n+s+i-1}$. Since $CK\times CL\cong C(K*L)$ and
$K*\sigma^0_0=CK$ for any polyhedra $K$ and $L$, it follows that
$$
Y^n\times I^{s+i}=(C\sigma^0_2)^n\times(C\sigma^0_0)^{s+i}\cong
\underbrace{C\dots C}_{s+i+1\text{ times}} (\sigma^0_2)^{*n}.
$$
If a polyhedron $K\not\emb S^d$ then the cone $CK\not\emb\R^{d+1}$
(because we work in piecewise linear category). So the non-embeddability in case
(1) follows from $(\sigma^0_2)^{*n}\not\emb S^{2n-2}$, which is proved in \cite{Kam32}
(or alternatively from $Y^n\not\emb S^{2n-1}$, which is proved in \cite{Sko07}
). \qed\enddemo

We are thus left with the proof of the non-embeddability in case
(2). To make it clearer, we anticipate it with considering
heuristically three simplest cases. Even a more visual way to express the main idea is given in \cite{Sko24}, where the so-called {\it linear} non-embeddability in the three examples is proved.

\smallskip

{\bf Example 1.} Let us first prove that the Kuratowski graph
$K_5$ not planar. Suppose to the contrary that $K_5\subset\R^2$.
Let $O$ be a vertex of $K_5$ and $D^2$ a small disc with the center
$O$. Then the intersection $K_5\cap\partial D^2$ consists of 4
points. Denote them by $A$, $B$, $C$, $D$, in the order along the
circle $\partial D^2$. Note that the pairs $A,C$ and $B,D$ are the
endpoints of two disjoint arcs contained in $K_5-\delet D^2$, and,
consequently, in $\R^2-\delet D^2$. Then the cycles $OAC, OBD\subset
K_5$ intersect each other ``transversely'' at the single point $O$,
which is impossible in the plane. So $K_5\not\emb\R^2$.
\smallskip

{\bf Example 2.} Now let us outline why $K_5\times
K_5\not\emb\R^4$. (Another proof is given in \cite{Um78}.) Recall
that if $K$ is a polyhedron with a fixed cell decomposition and $O\in K$ is a vertex, then {\it
the star} $\St O$ is the union of all closed cells of $K$
containing $O$, and {\it the link} $\Lk O$ is the union of all
cells of $\St O$ not containing $O$. In our previous example, $\Lk
O$ consists of 4 points and the proof is based on the fact that
there are two pairs of points of $\Lk O$ linked in $\partial D^2$.

Now take $K=K_5\times K_5$. Suppose to the contrary that $K\subset\R^4$. Let $O$ be a vertex of $K$ and $D^4$  be a small disc with the center $O$. Without loss of generality, the intersection $K\cap\partial D^4 =\Lk O\cong K_{4,4}$. So by the
Sachs theorem from the introduction, any embedding $\Lk O\emb \partial D^4$ has a
pair of linked cycles $\alpha,\beta\subset\Lk O$. Two linked cycles in $\partial D^4$ cannot bound two disjoint non-self-intersecting surfaces in $\Bbb{R}^4-\delet D^4$. If we construct two
such surfaces in $K-\St O$, then we get a contradiction, and
thus prove that $K\not\emb\R^4$. This construction is simple; see the proof of Lemma~2 below for details.

Analogously it can be shown that
$\sigma^2_6\not\emb\R^4$ (another proof is given in \cite{Kam32}.)
\smallskip

{\bf Example 3.} Let us show that $K_5\times S^1\not\emb\R^3$.
(Another proof was given by Tom Tucker; the simplest proof is analogous to Example~2 but we wish to illustrate another idea now.) Suppose that $K_5\times
S^1\emb\R^3$; then by (*) we have $K_5\times S^1\times S^1\emb\R^4$. But
$K_5\emb S^1\times S^1$, so $K_5\times K_5\emb\R^4$, which
contradicts Example~2.

\bigskip
{\bf Proof of the non-embeddability in case (2) modulo a lemma.}
Let $K$ and $L$ be two polyhedra, and a cell decomposition of $K$ is fixed. We say that a map $f:K\to L$ is {\it an almost embedding} \cite{FKT94}, if:

-  for any two {\it disjoint} closed cells $a,b\subset
K$ of the fixed decomposition we have $fa\cap fb=\emptyset$; and

- $f$ is piecewise linear on some subdivision of the fixed decomposition of $K$.

A nontrivial example of an almost embedding $K_5\to K_{3,3}$ is shown in Fig.~1. The following lemma is a ``half'' of the Menger conjecture up to replacement of ``almost embeddability'' by ``embeddability''. Hereafter the fixed cell decomposition is obvious and is not described explicitly.

\proclaim{Lemma~2} (For $n=2$, see \cite{Um78}) The polyhedron
$(K_5)^n$ is not almost embeddable into $\R^{2n}$.
\endproclaim

\demo{Proof of the non-embeddability in case $(2)$ of Theorem~1
modulo Lemma~2}
First we reduce the theorem to the particular case $s=0$ analogously to Example~3. Indeed, assume that $s>0$ and Theorem~1 does not hold in case (2) for a product $G\times (S^1)^s=G_1\times \dots\times G_n\times (S^1)^s$, i.e., $G\times (S^1)^s\hookrightarrow\R^{2n+s}$. Then by assertion~(*) from the proof of the embeddability in Theorem~1 it follows that
$$
G\times (K_5)^s\hookrightarrow
G\times (S^1)^{2s}      \hookrightarrow
\R^{2n+2s}.
$$
The composition of the two embeddings is an embedding of a product containing no factors homeomorphic to $S^1$. The existence of the latter embedding contradicts to the case $s=0$ of Theorem~1 because the theorem gives the dimension $d=2n+2s+1$ for the product $G\times (K_5)^s$.
The obtained contradiction reduces the theorem to the particular case $s=0$, which is considered now.

By the Kuratowsky graph planarity criterion any
nonplanar graph contains a subgraph homeomorhic to $K_5$ or
$K_{3,3}$. So we may assume that each $G_k$ is either $K_5$ or
$K_{3,3}$.
Now we are going to replace all the graphs $K_{3,3}$ by $K_5$-s.

Note that $K_5$ is almost embeddable to $K_{3,3}$ (Fig.~1).
Indeed, map a vertex of $K_5$ into the midpoint of an edge of
$K_{3,3}$ and map the remaining four vertices bijectively onto the four
vertices of $K_{3,3}$ not belonging to this edge. Then map each
edge $e$ of $K_5$ onto the shortest (with respect to the number of
vertices) arc in $K_{3,3}$, joining the images of the endpoints of $e$,
and the almost embedding is constructed.

A product of almost embeddings is again an almost
embedding, thus we get an almost embedding $(K_5)^n\to G_1\times\dots\times G_n$. Assume that there is an embedding $G_1\times\dots\times G_n\to\R^{2n}$. The composition of the almost embedding and the embedding is an almost embedding
$(K_5)^n\to \R^{2n}$, which contradicts to Lemma~2.
Thus the non-embeddability in case (2) of Theorem~1 follows from Lemma~2. \qed\enddemo

\bigskip
$
\epsfxsize=8.5cm\epsfbox{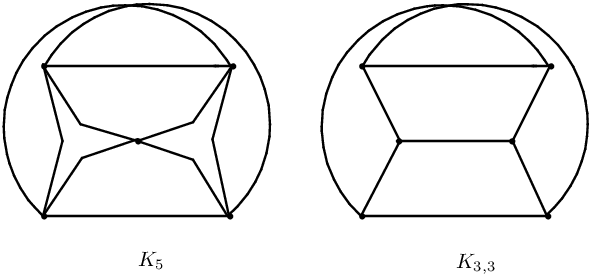}
$

\nopagebreak\smallskip
Fig. 1.
\smallskip

\bigskip
\bigskip


{\bf A few standard auxiliary facts.}
For the proof of Lemma 2, we need the following notions.
Let $A$ be a compact $n$-dimensional polyhedron. Its {\it boundary modulo $2$}, denoted $\partial A$,
is the union of all $(n-1)$-dimensional simplices contained in an odd number of $n$-dimensional ones, for some triangulation. 
Let $f:A\to\R^{2n}$, $g:B\to\R^{2n}$ be a pair of piecewise-linear maps such that $\dim A=\dim B=n\ne 0$ and $f\partial A\cap g B=f A\cap g\partial B=\emptyset$.
We say that the {\it intersection index} $fA\cap gB$ {\it is even}
(respectively, {\it odd}) and put $fA\cap gB:=0\in\Z/2\Z$ (respectively, $1\in\Z/2\Z$),
if the number of points in the set $\bar f A\cap \bar g B$ is even (respectively, odd),
for a general position pair of piecewise-linear maps $\bar f:A\to\R^{2n}$ and $\bar g:B\to\R^{2n}$ close to the pair of maps $f$ and $g$.

Let us explain the precise meaning of general position here. Take triangulations of $A$ and $B$ such that the maps $f:A\to\R^{2n}$, $g:B\to\R^{2n}$ are {\it linear}, i.e., linear on each simplex of the triangulations. Identify a pair of linear maps $\bar f:A\to\R^{2n}$, $\bar g:B\to\R^{2n}$ with a point $(\bar f,\bar g)\in\R^{2nN}$, where $N$ is the total number of vertices in the triangulations of $A$ and $B$. We say that $(\bar f,\bar g)$ is {\it $\varepsilon$-close} to $(f,g)$ if $|\bar fa-fa|<\varepsilon$ and $|\bar gb-gb|<\varepsilon$ for any vertices $a\in A$ and $b\in B$. We say that a property $P(\bar f,\bar g)$ holds {\it
for a general position pair of piecewise-linear maps $(\bar f,\bar g)$ close to $(f,g)$}, if there is $\varepsilon>0$ such that for each pair of triangulations of $A$ and $B$ such that $f$ and $g$ are linear, the property $P(\bar f,\bar g)$ holds for almost all $(\bar f,\bar g)\in\R^{2nN}$ that are $\varepsilon$-close to $(f,g)$.

We are going to use the following simple well-known result.

\proclaim{Parity Lemma 3}
If $A$ and $B$ are compact polyhedra of dimension $n\ne 0$ with
$\partial A=\partial B=\emptyset$, then for any 
piecewise-linear maps $f:A\to\R^{2n}$ and $g:B\to\R^{2n}$ the intersection index $fA\cap gB$ is even.
\endproclaim

The lemma 
just means vanishing of the intersection form in the homology of $\R^{2n}$ modulo 2. But it is simpler to prove it directly. See also \S4.5 
in~\cite{Sko24} 
and Remark 4.7.3c in~\cite{Sko23}.

\demo{Proof of Lemma~3} Take triangulations of $A$ and $B$ 
and identify pairs of linear maps $\bar f:A\to\R^{2n}$, $\bar g:B\to\R^{2n}$ with points in $\R^{2nN}$. 
For a pair of simplices $a\subset A$, $b\subset B$ and almost all $(\bar f,\bar g)\in\R^{2nN}$,
the set $\bar fa\cap \bar gb$ has at most one point, if $\dim a=\dim b=n$, and is empty, if $\dim a+\dim b=2n-1$. Indeed, otherwise 
$\bar fa$ and $\bar gb$ lie in one hyperplane in $\R^{2n}$, which is only possible for 
$(\bar f,\bar g)$ on certain algebraic hypersurface in $\R^{2nN}$.
Analogously, for almost all $(\bar f,\bar g)\in\R^{2nN}$, the set $\bar f A\cap \bar g B$ is in bijection
with the set $S$ of pairs of simplices $a\subset A,b\subset B$ such that $\bar f a\cap \bar g b\ne\emptyset $. 
In what follows we often omit the phrase ``for almost all $(\bar f,\bar g)\in\R^{2nN}$''.

It remains to prove that 
$S$ has an even number of elements.
Extend $\bar f:A\to \R^{2n}$ linearly to the cone $CA$ so that the cone vertex is mapped to the origin of $\R^{2n}$. The extension is still denoted by $\bar f$. Analogously to the above, for each pair of simplices $a\subset A$, $b\subset B$, 
the set
$\bar fCa\cap \bar gb$, if nonempty, is a single straight line segment when $\dim a=\dim b=n$, and a single point when $\dim a+\dim b=2n-1$; it is empty when $\dim a+\dim b<2n-1$.
Construct the following graph. The set of vertices is the set of pairs of simplices $a\subset CA,b\subset B$ such that $\dim a+\dim b=2n$ and $\bar f a\cap \bar g b\ne\emptyset $. Two distinct vertices $(a_1,b_1)$ and $(a_2,b_2)$ are joined by an edge, if there are top-dimensional simplices $a\subset CA$ and $b\subset B$ such that $a\supset a_1,a_2$ and $b\supset b_1,b_2$. In this case, $\bar f a\cap \bar g b$ is the segment joining $\bar f a_1\cap \bar g b_1$ and $\bar f a_2\cap \bar g b_2$; in particular, the latter points are distinct (otherwise,
either $\bar f\left|_a\right.$ is non-injective and $1=\dim \bar fa\cap \bar gb=\dim \bar f\partial a\cap \bar gb=0$, or $\bar g\left|_b\right.$ is non-injective, or $\bar f (a_1\cap a_2)\cap \bar g (b_1\cap b_2)\ne\emptyset$).
There is no other vertex $(a_3,b_3)$ with $a_3\subset a,b_3\subset b$ because the segment has just two endpoints. Thus the set of edges is in a bijection with the set of pairs of top-dimensional simplices $a\subset CA$ and $b\subset B$ such that $\bar f a\cap \bar g b\ne\emptyset$.
Hence the degree of a vertex $(a,b)$ is odd if and only if $a\subset A$, i.e., $(a,b)\in S$ (because $\partial (CA)=A\cup C\partial A$ and $\partial A=\partial B=\emptyset$). Thus $S$ has an even number of elements. 
\qed\enddemo


We need a well-known formula for the change of the intersection index upon a homotopy (Lemma~$3'$ below). 
Let $|S|$ be the number of elements in a set $S$.
For a map $F:A\times I\to\R^{2n}$, denote $F_t(x):=F(x,t)$. Let $F:A\times I\to\R^{2n}$, $G:B\times I\to\R^{2n}$ be piecewise-linear maps such that $\dim A=n$, $\dim B=n-1$ (or $B=\emptyset$), and $F_t\partial A\cap G_t B=F_t A\cap G_t\partial B=F_{\tau}A\cap G_{\tau} B=\emptyset$ for each $t\in I,\tau\in\{0,1\}$. We say that
the {\it intersection index} $\bigcup_{t\in I} F_t A\cap G_t B$ {\it of homotopies is even} (respectively, {\it odd}) and put $\bigcup_{t\in I} F_t A\cap G_t B:=0$ (respectively, $1$), 
if
the number $\left|\bigcup_{t\in I}\left(\bar F_t A\cap \bar G_t B\right)\right|$ is even (respectively, odd), for a general position pair of piecewise-linear maps $\bar F:A\times I\to\R^{2n}$ and $\bar G:B\times I\to\R^{2n}$ close to $F$ and $G$.
(This number equals the number of common points of the images of the maps $\bar F^{\uparrow}(x,t):=(\bar F(x,t),t)$ and $\bar G^{\uparrow}(y,t):=(\bar G(y,t),t)$ in $\R^{2n}\times I$.)

\proclaim{Lemma $3'$} Let $A$ and $B$ be compact polyhedra of dimension $n\ne 0$. Let $F:A\times I\to\R^{2n}$ and $G:B\times I\to\R^{2n}$ be piecewise-linear maps such that 
$F_t\partial A\cap G_t \partial B=F_\tau \partial A\cap G_\tau B=F_\tau A\cap G_\tau \partial B=\emptyset$ for each $t\in I,\tau\in\{0,1\}$. Then for a general position pair of piecewise-linear maps $(\bar F,\bar G)$ close to $(F,G)$,
we have
$$
|\bar F_1 A\cap \bar G_1 B|-|\bar F_0 A\cap \bar G_0 B|=
\left|\bigcup_{t\in I}(\bar F_t \partial A\cap \bar G_t B)\right| +\left|\bigcup_{t\in I}(\bar F_t A\cap \bar G_t \partial B)\right|
\mod 2.\qquad(3)
$$
Moreover, the intersection index is well-defined, i.e., is always either even or odd,
and therefore,
$$
F_1 A\cap G_1 B-F_0 A\cap G_0 B=
\bigcup_{t\in I} F_t \partial A \cap G_t B
+\bigcup_{t\in I} F_t A\cap G_t \partial B\mod 2.
$$
\endproclaim

\demo{Proof of Lemma $3'$}
Eq.~(3) is proved analogously to 
Lemma~3. Indeed, take triangulations of $A\times I$ and $B\times I$ and identify pairs of linear maps $\bar F:A\times I \to\R^{2n}$, $\bar G:B\times I\to\R^{2n}$ with points in $\R^{2nN}$. For a pair of simplices $a\subset A\times I$, $b\subset B\times I$ and almost all $(\bar F,\bar G)\in\R^{2nN}$,
the intersection $\bar F^{\uparrow}\delet a\cap \bar G^{\uparrow}\delet b$ is a single open interval of a straight line or empty when $\dim a=\dim b=n+1$,
a single point or empty when $\dim a+\dim b=2n+1$ or $\dim a+\dim b=2n$ with $a\subset A\times \{t\}$, $b\subset B\times \{t\}$ for some $t\in I$; and the intersection is empty otherwise.
For almost all $(\bar F,\bar G)\in\R^{2nN}$, we have $\left|\bigcup_{t\in I}(\bar F_t \partial A\cap \bar G_t B)\right|=\left|\bar F^{\uparrow} (\partial A\times I)\cap \bar G^{\uparrow} (B\times I)\right|$.
Take $(\bar F,\bar G)$ close enough to $(F,G)$ so that
$\bar F_t\partial A\cap \bar G_t \partial B=\bar F_\tau \partial A\cap \bar G_\tau B=\bar F_\tau A\cap \bar G_\tau \partial B=\emptyset$ for each $t\in I,\tau\in\{0,1\}$. Then Eq.~(3) is equivalent to the set $\left(\bar F^{\uparrow}\partial(A\times I)\cap \bar G^{\uparrow}(B\times I)
\right) \cup
\left(\bar F^{\uparrow}(A\times I)\cap \bar G^{\uparrow}\partial(B\times I)
\right)$ having an even number of elements.
We may assume that the latter set
is in bijection with the set $S$ of pairs of simplices $a\subset A\times I,b\subset B\times I$, where $a\subset \partial(A\times I)$ or $b\subset \partial(B\times I)$, such that $\bar F^{\uparrow}\delet a\cap \bar G^{\uparrow}\delet b\ne\emptyset$.
In what follows, we may impose assumptions, which are automatic for almost all $(\bar F,\bar G)\in\R^{2nN}$, without listing them explicitly.

To show that $|S|$ is even, construct the following graph. The set of vertices is the set of pairs of simplices $a\subset A\times I,b\subset B\times I$ such that $\bar F^{\uparrow} \delet a\cap \bar G^{\uparrow}\delet b$ is a single point.
Two distinct vertices $(a_1,b_1)$ and $(a_2,b_2)$ are joined by an edge, if there are top-dimensional simplices $a\subset A\times I$, $b\subset B\times I$ such that $a\supset a_1,a_2$ and $b\supset b_1,b_2$. In this case, $\bar F^{\uparrow} a\cap \bar G^{\uparrow} b$ is the segment joining $\bar F^{\uparrow} \delet a_1\cap \bar G^{\uparrow} \delet b_1$ and $\bar F^{\uparrow} \delet a_2\cap \bar G^{\uparrow} \delet b_2$ (in particular, the latter points are distinct because the images of the interiors of faces of $a$ are disjoint unless $\dim\bar F^{\uparrow} a\le n$).
Thus the set of edges is in a bijection with the set of pairs of top-dimensional simplices $a\subset A\times I$, $b\subset B\times I$ such that $\bar F^{\uparrow} \delet a\cap \bar G^{\uparrow} \delet b\ne\emptyset$.
Thus the set of vertices of odd degree is precisely $S$.
Thus $|S|$ is  even, and (3) holds.

Eq.~(3) implies that the intersection index $f A\cap g B$ is well-defined (if $f:A\to\R^{2n}$ and $g:B\to\R^{2n}$ satisfy $f \partial A\cap g B=f  A\cap g \partial B=\emptyset$). Indeed, let $F(x,t)=f(x)$ and $G(x,t)=g(x)$ be constant homotopies.
Take two triangulations of $A$ such that $f$ is linear on both and two maps $\bar F_0:A\to\R^{2n}$ and $\bar F_1:A\to\R^{2n}$ linear on the first and the second triangulation respectively. Extend the two triangulations to a triangulation of $A\times I$ such that $F$ is linear and the two maps to a linear map $\bar F:A\times I\to\R^{2n}$. 
Construct a triangulation of $B\times I$ and maps $G,\bar G:B\times I\to\R^{2n}$ analogously. All pairs  $(\bar F,\bar G)$ linear on our triangulations can be obtained by this construction.
If $(\bar F,\bar G)$ is close 
to $(F,G)$, then the right side of (3) vanishes, thus for some $\varepsilon>0$ the number
$|\bar F_0 A\cap \bar G_0 B|$ has the same parity for
all triangulations of $A$ and $B$ and almost all $(\bar F_0,\bar G_0)$, $\varepsilon$-close to $(f,g)$.

As a consequence, if linear maps $f:A\to\R^{2n}$ and $g:B\to\R^{2n}$ intersect 'transversely', i.e. each point of $f A\cap g B$ has a single $f$-preimage and a single $g$-preimage, both lying in the interior of top-dimensional simplices, then the intersection index $f A\cap g B=|f A\cap g B|\mod2$.

Finally, the intersection index of homotopies
$\bigcup_{t\in I} F_t A\cap G_t B$ 
is well-defined (if $F:A\times I\to\R^{2n}$ and $G:B\times I\to\R^{2n}$ satisfy the assumptions in its definition; in particular, now $\dim B=n-1$ or $B=\emptyset$) because 
$$\bigcup_{t\in I} F_t A\cap G_t B
=\left|\bigcup_{t\in I}\left(\bar F_t A\cap \bar G_t B\right)\right|
=\left|\bar F^{\uparrow} (A\times I)\cap \bar G^{\Uparrow} C\right|
=\bar F^{\uparrow} (A\times I)\cap \bar G^{\Uparrow} C
=F^{\uparrow} (A\times I)\cap G^{\Uparrow} C\mod2
$$
for a general position pair 
$(\bar F,\bar G)$ close to $(F,G)$.
Here the first equality is the definition. The second one is obvious, 
where we denote $C:=C(B\times I)\cong B\times I^2/B\times I\times\{0\}$ and
$\bar G^{\Uparrow}(y,t,s):=(2s\bar G^\uparrow(y,t),1-2s) 
\subset\Bbb{R}^{2n+1}\times\Bbb{R}$ for all $y\in B$ and $t,s\in I$.
The third one follows from the previous paragraph.
The last equality is obtained by applying~(3) to the rectilinear homotopy between $(\bar F^{\uparrow}, \bar G^{\Uparrow})$ and $(F^{\uparrow}, G^{\Uparrow})$.
\qed\enddemo

\bigskip

{\bf Completion of the proof of Theorem~1.} To complete the proof, it remains to prove Lemma~2. It
will be deduced from the following generalization of Lemma~1.

\proclaim{Lemma~$1'$} Let $L=(\sigma^0_3)^{*n}$. Then for any almost
embedding $f:CL\to\R^{2n}$ there exist two disjoint $(n-1)$-dimensional spheres
$\alpha,\beta\subset L$ such that the intersection index
$fC\alpha\cap fC\beta$ is odd.
\endproclaim

The cone $CL$ is considered here instead of $L$ itself. This auxiliary cone is essentially used only in the proof of {\it almost} non-embeddability in Lemma~2, not just non-embeddability.

\demo{Proof of Lemma~2 modulo Lemma~$1'$} Assume that there exists
an almost embedding $f:K=K_5\times\dots \times K_5\to\R^{2n}$. Let
$O=O_1\times\dots\times O_n$ be a vertex of $K$. By the well-known
formulae, we get
$$\Lk O\cong \Lk O_1*\dots*\Lk O_n
\text{ and } \St O=C\Lk O\cong C(\sigma^0_3)^{*n}.
$$

Let $\alpha,\beta\subset\Lk O$ be a pair of $(n-1)$-spheres given
by Lemma~$1'$. Identify $\Lk O$ and $\Lk O_1*\dots*\Lk O_n$.
Since $\alpha$ and $\beta$ are disjoint, it follows that for each
$k=1,\dots,n$ the sets $\alpha\cap\Lk O_k$ and $\beta\cap \Lk O_k$
are disjoint.
Each of the sets $\alpha\cap\Lk O_k$ and $\beta\cap \Lk O_k$ contains more than $1$ point because one of the spheres $\alpha$ and $\beta$ would be a cone otherwise. Thus each of the sets consists of exactly 2 points. By definition,
put $\{A_k,C_k\}:=\alpha\cap\Lk O_k$ and
$\{B_k,D_k\}:=\beta\cap\Lk O_k$. Consider two $n$-tori
$$
T_\alpha=O_1 A_1 C_1\times\dots\times O_n A_n C_n
\qquad\text{and}\qquad
T_\beta=O_1 B_1 D_1\times\dots\times O_n B_n D_n
$$
contained in $K$.

Clearly, $T_\alpha\supset C\alpha$, $T_\beta\supset C\beta$ and
$T_\alpha\cap T_\beta=O$. Since $f$ is an almost embedding, it
follows that the intersection index $fT_\alpha\cap fT_\beta=fC\alpha\cap fC\beta$. So $f
T_\alpha\cap f T_\beta$ is odd by the choice of $\alpha$ and $\beta$.
By Parity Lemma~3 we obtain a contradiction, hence $K$ is not almost embeddable into $\R^{2n}$.
\qed\enddemo

\demo{Proof of Lemma $1'$} The proof is similar to that of
Conway--Gordon--Sachs theorem and applies the idea of
\cite{Kam32}, 
with a slightly more refined obstruction. The reader
can restrict attention to the case when $n=2$ and obtain an
alternative proof of the Sachs theorem. (The proof for $n>2$ is
analogous to that for $n=2$.)

We show that for any $(n-1)$-simplex $c$ of $L$ and any almost
embedding $f:CL\to\R^{2n}$ there exist a pair of disjoint
$(n-1)$-spheres $\alpha,\beta\subset L$ such that $\alpha\supset
c$ and the intersection index $fC\alpha\cap fC\beta$ is odd.

For an almost embedding $f:CL\to\R^{2n}$, let $$v(f)=\sum
(fC\alpha\cap fC\beta)\mod2$$ be {\it the Van Kampen obstruction}
to linkless embeddability. Here the sum is over all pairs of
disjoint $(n-1)$-spheres $\alpha,\beta\subset L$ such that
$c\subset\alpha$. It suffices to prove that $v(f)=1$. Our proof is
in 2 steps: first we show that $v(f)$ does not depend on $f$;
then we calculate $v(f)$ for certain `standard' embedding $f:CL\to\R^{2n}$.

Let us prove that $v(f)$ does not depend on $f$ \cite{cf. Kam32,
CG83}. Take any two almost embeddings $F_0,F_1:CL\to\R^{2n}$. By
general position in piecewise linear category, there exists a piecewise linear homotopy $F:CL\times I\to\R^{2n}$ between them such that

1) there is a finite number of {\it singular moments} $t$, i.e., $t\in I$ such that $F_t$ is not an almost embedding;

2) for a singular $t$, there is a unique pair of disjoint
$(n-1)$-simplices $a,b\subset L$ such that $F_{t}Ca\cap
F_{t}b\ne\emptyset$;

3) the intersection $F_{t}Ca\cap F_{t}b$ is 'transversal in time', i.e.,
$F_{t}^{-1}(F_{t}Ca\cap F_{t}b)\times t$ consists of exactly two points, and $\left.F\right|_{Ca\times I}$ and $\left.F\right|_{b\times I}$ are smooth at those points.

Consider a singular moment $t$ and the pair $a,b$ of simplices given by condition 2). Conditions 3) and 1) imply that the intersection index of homotopies
$\bigcup_{\tau\in [t-\varepsilon,t+\varepsilon]} F_\tau Ca\cap F_\tau b$ is odd
for small enough $\varepsilon>0$.
Then by Lemma~$3'$, for a pair of
disjoint $(n-1)$-spheres $\alpha,\beta\subset L$,
the intersection index $F_tC\alpha\cap F_tC\beta$ changes with the
increasing of $t$ if and only if either $\alpha\supset a$,
$\beta\supset b$ or $\alpha\supset b$, $\beta\supset a$. Such
pairs $(\alpha,\beta)$ satisfying the condition $\alpha\supset c$
are called {\it critical}. If $c\cap (a\cup b)=\emptyset$, then
there are exactly $2$ critical pairs. Indeed, we have either
$\alpha\supset a\cup c$ or $\alpha\supset b\cup c$. Each of these
two conditions determines a unique critical pair. If $c\cap (a\cup
b)\ne\emptyset$, then there are two distinct vertices $v,w\in
L-(a\cup b\cup c)$ belonging to the same copy of $\sigma^0_3$.
Then there is an involution without fixed points on the set of
critical pairs. Indeed, $\Z/2\Z$ acts on the set of vertices of
$L$ by interchanging $v$ and $w$, and it also acts on the set of
critical pairs, because $v,w\notin a\cup b\cup c$. So the number
of critical pairs is 
even, thus 
$v(F_0)=v(F_1)$.

$
\epsfxsize=5cm\epsfbox{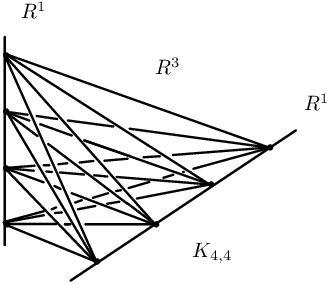}
$
\smallskip
Fig. 2.
\bigskip


Now let us prove that $v(f)=1$ for certain `standard' embedding
$f:CL\to\R^{2n}$ (Fig.~2). To define the standard embedding
$f:CL\to\R^{2n}$, take a general position collection of $n$ lines
in $\R^{2n-1}\subset\R^{2n}$. For each $k=1,\dots,n$ take a
quadruple $\sigma_k$ of distinct points on $k$-th line. Taking the
join of all $\sigma_k$, we obtain an embedding $L\to\R^{2n-1}$.
The standard embedding $f:CL\to\R^{2n}$ is defined to be the cone
over this embedding. Further we omit $f$ from the notation of
$f$-images.

Clearly, for a pair of disjoint $(n-1)$-spheres
$\alpha,\beta\subset L$, the intersection index $fC\alpha\cap
fC\beta$ has the same parity as the linking number $\lk(\alpha,\beta)$ in $\R^{2n-1}$.
Indeed, let $H$ be the half-space bounded by $\R^{2n-1}$ and containing the cone vertex.
Take a pair of triangulations of $C\alpha$ and $C\beta$ and
a pair of linear embeddings $(\bar f\colon C\alpha\to H, \bar g\colon C\beta\to H)$ close enough to $(\left.f\right|_{C\alpha},\left.f\right|_{C\beta})$
such that $\bar f\alpha=\alpha$, $\bar g\beta=\beta$, and $\bar fa\cap\bar gb=\emptyset$ for each pair of simplices $a\subset C\alpha,b\subset C\beta$ unless $\dim a=\dim b=n$.
Applying Lemma~$3'$ to the rectilinear homotopy between $(\bar f,\bar g)$ and $(\left.f\right|_{C\alpha},\left.f\right|_{C\beta})$, we get $fC\alpha\cap
fC\beta=\bar fC\alpha\cap\bar gC\beta=|\bar fC\alpha\cap\bar gC\beta|\mod 2$. 
Passing to a subdivision, if necessary, we guarantee that for each simplex $b\subset C\beta$ there is at most one simplex $a\subset C\alpha$ such that $\bar fa\cap\bar gb\ne\emptyset$. Then $\beta$ is homologous in $H-\bar fC\alpha$ to the sum of the boundaries $\partial b$ of all simplices $b\subset C\beta$ such that $\bar fa\cap\bar gb\ne\emptyset$ for some simplex $a\subset C\alpha$. A simple computation shows that each such $\partial b$ contributes $1$ to $\lk(\alpha,\beta)\mod2$, hence 
$|\bar fC\alpha\cap\bar gC\beta|=\lk(\alpha,\beta)\mod2$.

Let us show that
$\lk(\alpha,\beta)=1\mod2$ if and only if for each $k=1,\dots,n$
the 0-spheres $\alpha\cap \sigma_k$ and $\beta\cap \sigma_k$ are
linked in $k$-th line.
Indeed, if, say, the two points $\alpha\cap \sigma_k$ lie between the two points of $\beta\cap\sigma_k$ for some $k$, then take the segment $I_\alpha$ with $\partial I_\alpha=\alpha\cap \sigma_k$. The piecewise linear disk
$$D_\alpha=(\alpha\cap \sigma_1)*\dots*(\alpha\cap \sigma_{k-1})*I_\alpha*
(\alpha\cap \sigma_{k+1})*\dots*(\alpha\cap \sigma_n)$$ spans $\alpha$ and is disjoint with $\beta$, hence $\lk(\alpha,\beta)=0\mod2$. If $\alpha\cap \sigma_k$ and $\beta\cap \sigma_k$ are linked for each $k=1,\dots,n$, then take the minimal segment $I_k$ containing $\sigma_k$.
The complement to
$\alpha\cap \sigma_k$ in $I_k$ deformationally retracts to $\beta\cap \sigma_k$, thus the complement to $\alpha$ in the simplex $I_1*\dots *I_n$ deformationally retracts to $\beta$.
Then for any $(n-1)$-sphere $\gamma$ in the latter complement,
$\lk(\alpha,\gamma)$ is a multiple of $\lk(\alpha,\beta)$. Since there exists $\gamma$ with $\lk(\alpha,\gamma)=1\mod2$, it follows that $\lk(\alpha,\beta)=1\mod2$. (Another explanation is that $D_\alpha$ and $\beta$ intersect transversely at a single point, but we avoid piecewise linear transversality.)

Now it is obvious that there exists a unique pair
$\alpha,\beta$ such that $\alpha\supset c$ and $fC\alpha\cap
fC\beta=1\mod2$. So $v(f)=1$, which proves the lemma. \qed
\enddemo

We conclude the paper by the proof of Theorem~1 in topological category
(due to the referee):

\demo{Proof of Theorem~1 in topological category} For codimension
$\ge3$, the assertion of Theorem~1 in topological category follows from
the one in piecewise linear one, because by a theorem of Bryant \cite{Bry72}
each topological embedding of a polyhedron into a piecewise-linear manifold in codimension $\ge 3$ can be approximated by a piecewise linear embedding.

The cases of codimension 1 and 2 are reduced to codimension $\ge 3$ case analogously to Example 3. Indeed, assume that Theorem~1 does not hold for a polyhedron $K=G_1\times \dots G_n\times I^i\times (S^1)^s$. This means that $K\hookrightarrow_{TOP}\R^{d(K)-1}$, where 
$d=d(K)$ is the dimension
given by (1)--(2). 
By assertion~(*), 
$$
K\times K_5\times K_5\hookrightarrow_{TOP}
K\times (S^1)^4      \hookrightarrow_{TOP}
\R^{d(K)+3}.
$$
The composition of the two embeddings is now a codimension $\ge 3$ embedding. The existence of the latter contradicts to the codimension $\ge 3$ case of Theorem~1 in topological category because $d(K\times K_5\times K_5)=d(K)+4$.
The obtained contradiction proves the theorem.
 \qed\enddemo

\smallskip
{\bf Acknowledgements.} The author is grateful to Arkady Skopenkov
for permanent interest to this work and to the referee for useful
suggestions and a remark proving one of the author's conjectures.
The author thanks Emil Alkin for a discussion that inspired the addition of proofs on p.~4--5 to the updated version.

\enddocument
\end